\documentclass[11pt]{article}
\usepackage[a4paper]{anysize}\marginsize{2cm}{2cm}{2cm}{2cm}
\pdfpagewidth=\paperwidth \pdfpageheight=\paperheight
\usepackage{amsfonts,amssymb,amsthm,amsmath,eucal}
\usepackage{bbm}
\usepackage{tikz} 
\usepackage{mathrsfs}

\theoremstyle{plain}
\newtheorem{thm}{Theorem}[section]
\newtheorem{theorem}[thm]{Theorem}

\newtheorem{lemma}[thm]{Lemma}
\newtheorem{proposition}[thm]{Proposition}

\theoremstyle{definition}

\newtheorem{question}[thm]{Question}
\newtheorem{problem}[thm]{Problem}

\newtheorem{thevarthm}[thm]{\varthmname}

\newenvironment{varthm*}[1]{\trivlist\item[]{\bf #1.}\it}{\endtrivlist}

\renewcommand\geq{\geqslant}

\renewcommand\leq{\leqslant}

\newcommand\be{\begin{eqnarray*}}
\newcommand\ee{\end{eqnarray*}}

\newcommand\C{\mathbb C}

\newcommand\m{\mathfrak{m}}
\newcommand\p{\mathfrak{p}}
\renewcommand\P{\mathbb P}

\newcommand\newop[2]{\def#1{\mathop{\rm #2}\nolimits}}
\newop\log{log}
\newop\ord{ord}
\newop\Gal{Gal}
\newop\SL{SL}
\newop\Bl{Bl}
\newop\mult{mult}
\newop\imult{imult}
\newop\mass{mass}
\newop\Ass{Ass}
\newop\div{div}
\newop\codim{codim}
\newop\sing{sing}
\newop\Zeroes{Zeroes}

\newcommand\wtilde[1]{\widetilde{#1}}

\def\keywordname{{\bfseries Keywords}}%
\def\keywords#1{\par\addvspace\medskipamount{\rightskip=0pt plus1cm
\def\and{\ifhmode\unskip\nobreak\fi\ $\cdot$
}\noindent\keywordname\enspace\ignorespaces#1\par}}
\def\subclassname{{\bfseries Mathematics Subject Classification
(2000)}\enspace}
\def\subclass#1{\par\addvspace\medskipamount{\rightskip=0pt plus1cm
\def\and{\ifhmode\unskip\nobreak\fi\ $\cdot$
}\noindent\subclassname\ignorespaces#1\par}}

\def\endproof{\hspace*{\fill}\endproofsymbol\endtrivlist}

\def\endproofsymbol{\frame{\rule[0pt]{0pt}{6pt}\rule[0pt]{6pt}{0pt}}}

\begin{document}

\title{On codimension two flats in Fermat-type arrangements}
\author{Grzegorz Malara and Justyna Szpond}
\maketitle

\abstract{In the present note we study certain arrangements of codimension $2$ flats in projective
spaces, we call them \emph{Fermat arrangements}. We describe algebraic properties of their defining
ideals. In particular, we show that they provide counterexamples to an expected containment relation
between ordinary and symbolic powers of homogeneous ideals.}

\section{Introduction}
\label{SSsec:1}
   A Fermat-type arrangement of degree $n\geq 1$ of hyperplanes in projective space $\P^N$
   is given by linear factors of the polynomial
   $$F_{N,n}=F_{N,n}(x_0,\ldots,x_N)=\prod\limits_{0\leq i<j\leq N}(x_i^n-x_j^n).$$
   These arrangements sometimes appear under the name of Ceva arrangements in the literature,
   see e.g. \cite[Section 2.3.I]{BHH}. The name Fermat arrangement has been used
   for lines in $\P^2$ e.g. by Urzua, see \cite[Example II.6]{Urz08}.
   Fermat arrangements of lines have attracted recently considerable attention, see e.g. \cite{NagSec16},
   because of their appearance on the border line of the following fundamental problem.
   \begin{problem}[Containment problem]\label{pro:containment}
      Let $I$ be a homogeneous ideal in the polynomial ring $\C[x_0,\ldots,x_N]$.
      Determine all pairs of integers $m$ and $r$ such that the
      containment
      \begin{equation}\label{eq:containment}
         I^{(m)}\subset I^r
      \end{equation}
      between the symbolic and ordinary powers of the ideal $I$ holds.
   \end{problem}
   We recall that for $m\geq 0$ the $m$-th symbolic power of $I$ is defined as
   \begin{equation}\label{eq:symbolic power}
      I^{(m)}=\bigcap_{P\in\Ass(I)}\left(I^mR_P\cap R\right),
   \end{equation}
   where $\Ass(I)$ is the set of associated primes of $I$.
   A ground breaking result of Ein, Lazarsfeld and Smith
   \cite{ELS01} (rendered by  Hochster and Huneke
   in positive characteristic \cite{HoHu02})
   asserts that there is always
   the containment in \eqref{eq:containment} for $m$ and $r$ subject to the inequality
   \begin{equation}\label{eq:ELS HH}
      m\geq hr,
   \end{equation}
   where $h$ is the maximum of heights of all associated primes of $I$.
   The natural question: To which extend the bound in \eqref{eq:ELS HH}
   is sharp has fueled a lot of research in the last 15 years. Considerable
   attention has been given to the following question of Huneke.
\begin{question}[Huneke]\label{que:Huneke}
   Let $Z$ be a finite set of points in $\P^2$ and let $I$ be the homogeneous
   ideal defining $Z$. Is then
   \begin{equation}\label{eq:I3 in I2}
      I^{(3)}\subset I^2?
   \end{equation}
\end{question}
   Note that the containment $I^{(4)}\subset I^2$ follows in this situation
   directly from \eqref{eq:ELS HH}. On the other hand it is very
   easy to find sets of points in $\P^2$ for which
   the containment $I^{(2)}\subset I^2$ fails. Question \ref{que:Huneke} remained
   open for quite a long time. It is by now known that there are sets of
   points in $\P^2$ for which the containment in \eqref{eq:I3 in I2}
   fails. The first counterexample has been given by Dumnicki, Szemberg
   and Tutaj-Gasi\'nska in \cite{DST13}. This counterexample is provided
   by the set of all $12$ intersection points of the Fermat arrangement of $9$
   lines in $\P^2$ (i.e. $n=3$ in this arrangement). The paper \cite{DST13}
   suggested that arrangements with arbitrary $n\geq 3$ should provide further
   counterexamples. This has been worked out and verified to hold by Harbourne
   and Seceleanu, see \cite[Proposition 2.1]{HarSec15}. Whereas Fermat configurations
   of lines allow no deformations, a series of counterexamples allowing parameters
   has been presented recently in \cite{LBS17}. Apart of those series there are some
   sporadic counterexamples to Question \ref{que:Huneke} available. The nature
   of all these examples has not been yet fully understood.

   The statement of \eqref{eq:containment} does not restrict to ideals
   supported on points. In particular, Huneke's question
   can be reformulated for codimension $2$ subvarieties in the projective
   space of arbitrary dimension $N$.
\begin{question}[Huneke-type]
   Let $V$ be codimension $2$ subvariety in $\P^N$ and let $I$ be the homogeneous
   ideal defining $V$. Is then
   \begin{equation}\label{eq:I3 in I2 H-type}
      I^{(3)}\subset I^2?
   \end{equation}
\end{question}
   This question has been answered to the negative in \cite{MalSzp17}.
   More precisely, we showed \cite[Theorem 4.3]{MalSzp17}
   that the containment \eqref{eq:I3 in I2 H-type}
   fails for $V$ consisting of all lines in $\P^3$
   contained in at least $3$ hyperplanes among those defined by linear factors of $F_{3,n}$
   for all $n\geq 3$. In the present note, we investigate the ideals
   defining codimension $2$ linear flats of multiplicity at least $3$ cut out
   by the Fermat-type arrangement given by $F_{N,n}$ with $n\geq 3$, as well
   as ordinary and symbolic powers of these ideals.

\section{Notation and basic properties}
\label{SSsec:2}
   The bookkeeping of all data is quite essential for what follows.
   In the present section we establish the notation and prove some
   basics facts.

   By $x_0,x_1,\ldots,x_N$ we denote the coordinates in the projective space $\P^N$.
   We fix an integer $n\geq 3$. This integer is not present in the short hand notation
   introduced below. We hope that it will not lead to any confusion since we
   work always with $n$ fixed. We introduce the following bracket symbol.
   For an integer $1\leq k\leq N$ let $i_0,\ldots,i_k$ be $(k+1)$ mutually distinct
   elements in the set $\left\{0,\ldots,N\right\}$.
   \begin{equation}\label{eq:bracket def}
   [x_{i_1}\ldots x_{i_k}]:=\prod_{p<q}(x_{i_p}^n-x_{i_q}^n).
   \end{equation}
   Thus, in particular,
   $$F_{N,n}=[x_0\ldots x_N].$$
   The notation in \eqref{eq:bracket def} fulfills the antisymmetry condition.
   More precisely, we have for any pair $p,q$ such that $1\leq p<q\leq k$:
   \begin{equation}\label{eq:antisymmetry}
      [x_{i_1}\ldots x_{i_p}\ldots x_{i_q}\ldots x_{i_k}]=(-1)^{q-p}[x_{i_1}\ldots x_{i_q}\ldots x_{i_p}\ldots x_{i_k}].
   \end{equation}
\begin{lemma}[Expansion rule]\label{rul:expansion}
   For $k\geq 2$ there is
   $$[x_{i_0}\ldots x_{i_{k-1}}x_{i_k}]=[x_{i_0}\ldots x_{i_{k-1}}]\prod_{j=0}^{k-1}[x_{i_j}x_{i_k}].$$
\end{lemma}
\proof
   This follows straightforward from the definition in \eqref{eq:bracket def}.
\endproof
   For example
   $$[xyzw]=[xyz](x^n-w^n)(y^n-w^n)(z^n-w^n).$$
   We have also the following Laplace-type rule.
\begin{lemma}[Laplace expansion]\label{rul:Laplace}
   We have
   $$[x_{i_0}\ldots x_{i_k}]=\sum_{j=0}^k (-1)^{j+k} x_{i_0}^n\ldots \widehat{x_{i_j}^n}\ldots x_{i_k}^n[x_{i_0}\ldots \widehat{x_{i_j}}\ldots x_{i_k}].$$
   As usual $\widehat{a}$ means that the term $a$ is omitted.
\end{lemma}
\proof
In order to alleviate notation, we drop the double index notation. It's sufficient to show
   \begin{equation}\label{eq:Laplace}
      [x_{0}\ldots x_{k}]=\sum_{j=0}^k (-1)^{j+k}x_{0}^n\ldots\widehat{x_j^n}\ldots x_{k}^n[x_{0}\ldots\widehat{x_j}\ldots x_{k}].
   \end{equation}
    Both sides are polynomials of degree $\frac{k(k+1)}{2}n$. It is enough to show that the right hand side vanishes along all hyperplanes
    of the form $x_i=\delta x_j$, where $\delta$ is some root of $1$ of degree n.
    By symmetry it is enough to check this for $x_0=\delta x_1$. Then the right hand side in (\ref{eq:Laplace}) is
    $$(-1)^{k}x_{1}^n\ldots x_{k}^n[x_{1}\ldots x_{k}]+ (-1)^{k+1}x_{0}^n x_{2}^n\ldots x_{k}^n[x_{0}x_{2}\ldots x_{k}],$$
    since $x_0=\delta x_1$ we get $0$.

    Thus the right hand side of (\ref{eq:Laplace}) is equal to $\lambda \cdot [x_{0}\ldots x_{k}]$ for some $\lambda \in \C$.
    In order to establish $\lambda$, we evaluate at $x_0=0$, which gives
    $$\lambda(-x_{1}^n)\ldots (-x_{k}^n)[x_{1}\ldots x_{k}]=(-1)^{k}x_{1}^n\ldots x_{k}^n[x_{1}\ldots x_{k}],$$
    hence $\lambda=1$.
\endproof
   Also the next fact is very useful.
\begin{lemma}[Substitution rule]\label{rul:substitution}
   For any $u\in\left\{0,\ldots,N\right\}$ and $1\leq k\leq N$ there is
   $$[x_{i_0}\ldots x_{i_k}]=\sum_{j=0}^k[x_{i_0}\ldots x_{i_{j-1}}x_u x_{i_{j+1}}\ldots x_{i_k}].$$
\end{lemma}
   For example
   $$[xyz]=[wyz]+[xwz]+[xyw].$$
\proof
   In order to alleviate notation we drop the double index notation. It is clear that the statement
   is invariant under the symmetry group on the $(N+1)$ variables. Also it convenient to use \eqref{eq:antisymmetry}
   and write the assertion in the following form
   \begin{equation}\label{eq:substitution asymmetric}
      [x_{0}\ldots x_{k}]=\sum_{j=0}^k (-1)^{j}[x_u x_{0}\ldots x_{j-1}\widehat{x_j} x_{j+1}\ldots x_{k}].
   \end{equation}
   The argumentation is similar to that in proof of Lemma \ref{rul:Laplace}. Both sides in (\ref{eq:substitution asymmetric}) are
   polynomials of degree $\frac{k(k+1)}{2}n$. We substitute $x_0=\delta x_1$. Then the right hand side is
   $$[x_u x_1\ldots x_k]-[x_u x_0 x_2 \ldots x_k] $$
   which is clearly $0$. Thus we have
   $$\sum_{j=0}^k[x_{i_0}\ldots x_{i_{j-1}}x_u x_{i_{j+1}}\ldots x_{i_k}]=\lambda \cdot [x_{0}\ldots x_{k}],$$
   for some $\lambda \in \C$. In order to determine $\lambda$, we substitute $x_u=x_0$. Then
   $$[x_{0}\ldots x_{k}]=\lambda \cdot [x_{0}\ldots x_{k}],$$
   which implies $\lambda=1$.
\endproof
We conclude these preparations by another useful rule.
\begin{lemma}[Useful rule]\label{rul:useful rule}
   For $k\geq2$ and auxiliary variables $y_1,\ldots,y_k$ we have
   $$[x_{0}\ldots x_{k}]=\sum_{j=0}^k (-1)^j [x_{0}\ldots\widehat{x_j}\ldots x_{k}][x_j y_1]\ldots[x_j y_k] .$$
\end{lemma}
\proof
The proof parrots that of Lemma \ref{rul:Laplace} and Lemma \ref{rul:substitution} and is left to the reader.
In order to determine the constant $\lambda$ one might substitute $y_i=x_i$ for $i=1, \ldots, k$.
\endproof

\section{Fermat arrangements of codimension two flats}
\label{SSsec:3}
   In this section we study, for $n\geq 3$, the union $V_{N,n}$ of codimension $2$ flats $W$ in $\P^N$
   such that there are at least $3$ hyperplanes among those defined by the linear factors of $F_{N,n}$
   vanishing along $W$. Let $I_{N,n}$ be the radical ideal defining $V_{N,n}$. The set $V_{N,n}$
   is the union of $N+1$ cones with vertices in the coordinate points $E_i=(0:\ldots:\underbrace{1}_{\text{$i$}}\ldots:0)$
   over the sets $V_{N-1,n}(i)$ defined
   in the hyperplanes $H_i=\left\{\;x_i=0\;\right\}$. Let $I_{N-1,n}(i)$ be the ideal defining $V_{N-1,n}(i)$
   in the variables $x_0,\ldots,\widehat{x_i},\ldots,x_N$. The geometry
   of the arrangement implies the following relation between the defined ideals.
\begin{lemma}\label{lem:I_N,n as intersection}
   Keeping the notation above, we have for all $N\geq 3$
   $$I_{N,n}=\bigcap_{i=0}^N I_{N-1,n}(i).$$
\end{lemma}
   For the proof of the main Theorem \ref{thm:non-containment} we need a more direct
   description of ideals $I_{N,n}$ in terms of generators.
\begin{proposition}\label{prop:generators}
   Consider the ideal $I_{N,n}$ for some integers $N\geq 2$ and $n\geq 3$.\\
   a) Let $N=2M$ be an even number. Let $A=\left\{i_1,\ldots,i_M\right\}$ be a subset of $M$ elements
   in the set $\left\{0,1,\ldots,N\right\}$ and let $B=\left\{j_0,\ldots,j_{M}\right\}$
   be the complimentary set. The ideal $I_{N,n}$ is generated by all polynomials
   of the form
   $$g_A=x_{i_1}\ldots x_{i_M}[x_{i_1}\ldots x_{i_M}][x_{j_0}\ldots x_{j_{M}}].$$
   b) Let $N=2M+1$ be an odd number. Let $A=\left\{i_0,\ldots,i_M\right\}$ be a subset of $M+1$ elements
   in the set $\left\{0,1,\ldots,N\right\}$ and let $B=\left\{j_0,\ldots,j_{M}\right\}$
   be the complimentary set. The ideal $I_{N,n}$ is generated by all polynomials
   of the form
   $$g_A=x_{i_0}\ldots x_{i_M}[x_{i_0}\ldots x_{i_M}][x_{j_0}\ldots x_{j_{M}}].$$
\end{proposition}
\proof
   The proof goes by induction on $N$. The first step, $N=2$ has been shown in \cite[Lemma 2.1]{DST13}.
   Using the presentation in Lemma \ref{lem:I_N,n as intersection}, we will show that
   generators $g_A$ are contained in each of the intersecting ideals. To this end we study first
   \textbf{the case $N$ is even} with $N=2M$.

   Since everything is invariant under the permutation group, it suffices to work with the set
   $A=\left\{0,1,\ldots,M-1\right\}$. Then
   $$g_A=x_0x_1\ldots x_{M-1}[x_0\ldots x_{M-1}][x_M\ldots x_{2M}].$$

   We have the following two cases.
   Assume that $M\leq i\leq 2M$. Then the ideal $I_{N-1,n}(i)$ contains
   as a generator
   $$h_{A}=x_0\ldots x_{M-1}[x_0\ldots x_{M-1}][x_M\ldots\widehat{x_i}\ldots x_{2M}].$$
   It is easy to see that $g_A$ is divisible by $h_{A}$, indeed
   $$g_A=\pm\prod_{j=M}^{2M}[x_jx_i] h_{A}.$$
   Assume now that $0\leq i\leq M-1$. After renumbering the variables we can in fact assume
   that $i=0$. Let $A_j=\left\{1,\ldots,M-1,M+j\right\}$ for $j=0,\ldots,M$ and
   let $h_{A_j}$ be the corresponding generators of $I_{N-1,n}(i)$, i.e.
   $$h_{A_j}=x_1\ldots x_{M-1}x_{M+j}[x_1\ldots x_{M-1}x_{M+j}][x_M\ldots\widehat{x_{M+j}}\ldots x_{2M}].$$
   Then
   $$g_A=\sum_{j=0}^M(-1)^{j+M-1}x_0x_{M+j}^{n-1}[x_0x_1]\ldots[x_0x_{M-1}]h_{A_j}.$$
   To see this we will alter the right hand side of the above equality. First note that by Lemma \ref{rul:expansion} we have
   $$\sum_{j=0}^M(-1)^{j+M-1}x_0x_{M+j}^{n-1}[x_0x_1]\ldots[x_0x_{M-1}]x_1\ldots x_{M-1}x_{M+j}[x_1\ldots x_{M-1}x_{M+j}][x_M\ldots\widehat{x_{M+j}}\ldots x_{2M}]=$$
   $$x_0\ldots x_{M-1}\sum_{j=0}^M(-1)^{j+M-1}x_{M+j}^{n}[x_0x_1]\ldots[x_0x_{M-1}][x_1\ldots x_{M-1}][x_1x_{M+j}]\ldots[x_{M-1}x_{M+j}][x_M\ldots\widehat{x_{M+j}}\ldots x_{2M}].$$
   Again by the Expansion rule it reduces to:
   $$x_0\ldots x_{M-1}[x_0\ldots x_{M-1}]\sum_{j=0}^M(-1)^{j}[x_{M+j}0][x_{M+j}x_1]\ldots[x_{M+j}x_{M-1}][x_M\ldots\widehat{x_{M+j}}\ldots x_{2M}].$$
   By Lemma \ref{rul:useful rule} with $y_1=x_1$, $\ldots$, $y_{M-1}=x_{M-1}$, $y_M=0$ this expression reduces to
   $$x_0\ldots x_{M-1}[x_0\ldots x_{M-1}][x_M\ldots x_{2M}]=g_A.$$
   Now we pass to \textbf{the case $N$ is odd} with $N=2M+1$.

   Let $A=\left\{0,\ldots,M\right\}$ and
   $$g_A=x_0x_1\ldots x_{M}[x_0\ldots x_{M}][x_{M+1}\ldots x_{2M+1}].$$
   There are again two subcases. Assume that $0\leq i\leq M$. Then the ideal $I_{N-1,n}(i)$
   contains the generator
   $$g_{A'}=x_0x_1\ldots\widehat{x_i}\ldots x_{M}[x_0\ldots\widehat{x_i}\ldots x_{M}][x_{M+1}\ldots x_{2M+1}]$$
   with $A'=A\setminus\left\{i\right\}$.
   Then $g_A$ is divisible by $g_{A'}$, indeed
   $$g_A=\pm x_i[x_0x_i]\ldots\widehat{[x_ix_i]}\ldots[x_Mx_i] g_{A'}.$$
   For $M+1\leq i\leq 2M+1$ it suffices, up to renumbering the variables to consider $i=2M+1$.
   In the ideal $I_{N-1,n}(i)$ there are generators
   $$g_{A_j}=x_0x_1\ldots\widehat{x_j}\ldots x_{M}[x_0\ldots\widehat{x_j}\ldots x_{M}][x_jx_{M+1}\ldots x_{2M}]$$
   for $A_j=\left\{0,1,\ldots\widehat{j}\ldots,M\right\}$.
   Then
   $$g_A=\sum_{j=0}^M (-1)^jx_j[x_{M+1}x_{2M+1}]\ldots[x_{2M}x_{2M+1}]\cdot g_{A_j}.$$
   Again, we reduce the right hand side of this equality. To begin with we have
   $$\sum_{j=0}^M (-1)^jx_j[x_{M+1}x_{2M+1}]\ldots[x_{2M}x_{2M+1}]\cdot x_0x_1\ldots\widehat{x_j}\ldots x_{M}[x_0\ldots\widehat{x_j}\ldots x_{M}][x_jx_{M+1}\ldots x_{2M}]=$$
   $$x_0\ldots x_M\sum_{j=0}^M (-1)^j [x_0\ldots\widehat{x_j}\ldots x_{M}][x_jx_{M+1}]\ldots[x_jx_{2M}][x_{M+1}\ldots x_{2M}][x_{M+1}x_{2M+1}]\ldots[x_{2M}x_{2M+1}].$$
   Combining this with Lemma \ref{rul:expansion} and Lemma \ref{rul:Laplace} we get
   $$x_0\ldots x_M [x_0\ldots x_{M}][x_{M+1}\ldots x_{2M+1}]=g_A.$$
   Thus we have shown that in both cases every generator
   $$g_A\in \bigcap_{i=0}^N I_{N,n-1}(i).$$
   It remains to check that the ideal generated by all $g_A$ is indeed the whole ideal $I_{N,n}$. We leave this to a motivated reader.
\endproof
\section{The non-containment result}
   In this section we prove our main result.
\begin{theorem}\label{thm:non-containment}
   For arbitrary $N\geq 2$ and $n\geq 3$ there is
   $$I_{N,n}^{(3)}\not\subset I_{N,n}^2.$$
\end{theorem}
\proof
   It is convenient to abbreviate $I=I_{N,n}$.
   The polynomial $f:=f_{N,n}=[x_0\ldots x_{2M}]$ is contained in $I^{(3)}$
   by the Zariski-Nagata Theorem, see \cite[Theorem 3.14]{Eisenbud} for prime ideals
   and \cite[Corollary 2.9]{SidSul09} for radical ideals.
   Let $G$ denote the set of generators of the ideal $I$.

   The proof that it is not contained in $I^2$
   depends on the parity of the dimension $N$ of the ambient space.

   We handle first \textbf{the case $N=2M$}.
   Assume to the contrary that $f\in I^2$. Then there are
   polynomials $h_{g,g'}$ such that
   \begin{equation}\label{eq:1}
      f=\sum\limits_{g,g'\in G} h_{g,g'} gg'.
   \end{equation}
   Taking \eqref{eq:1} modulo $(x_0)$ we have
   \begin{equation}\label{eq:1 mod x0}
      \wtilde{f}=\sum\limits_{g,g'\in G} \wtilde{h_{g,g'}}\cdot \wtilde{g}\cdot \wtilde{g'},
   \end{equation}
   where $\wtilde{q}$ denotes the residue class of $q\in\C[x_0,\ldots,x_N]$ modulo $(x_0)$.
   Then
   $$\wtilde{f}=x_1^n\ldots x_{2M}^n[x_1\ldots x_{2M}].$$
   We focus now on the coefficient at the monomial
   $$\m=x_1^{2Mn}x_2^{(2M-1)n}\ldots x_{2M-1}^{2n} x_{2M}^n$$
   on both sides of equation \eqref{eq:1 mod x0}. This coefficient is $1$
   on the left hand side of \eqref{eq:1 mod x0}. It is easy to see that
   there is exactly one way to get this monomial expanding the product defining $\wtilde{f}$.

   Let $g\in G$ be a generator of $I$. By Proposition \ref{prop:generators} $g$ has the form
   $$g=x_{i_1}\ldots x_{i_M}[x_{i_1}\ldots x_{i_M}][x_{j_0}\ldots x_{j_M}],$$
   with all indices $i_1,\ldots,i_M,j_0,\ldots,j_M$ mutually distinct.
   If $0\in\left\{i_1,\ldots,i_M\right\}$, then the residue class of $g$ is zero.
   If it is in the second group of indices, then the residue class, after possible
   renumbering of indices, has the form
   \begin{equation}\label{eq:g wtilde}
      \wtilde{g}=x_{i_1}\ldots x_{i_M}x_{j_1}^n\ldots x_{j_M}^n[x_{i_1}\ldots x_{i_M}][x_{j_1}\ldots x_{j_M}].
   \end{equation}
   Note that we suppress the notation and write $x_i$ rather than $\wtilde{x_i}$.

   We will now analyze how the monomial $\m$ appears on the right hand side of \eqref{eq:1 mod x0}.
   To this end we run the following procedure starting with the variables with least powers
   in $\m$.
   \begin{itemize}
   \item The variable $x_{2M}$ has to be among the variables appearing with power $1$ in the
   product defining $\wtilde{g}$ and $\wtilde{g'}$ (variables indexed by the letter $i$) because its total power in $\m$ is restricted
   by $n$ and this is the only possibility to fulfill this condition.
   \item The variable $x_{2M-1}$ cannot then appear with power $1$ neither in $\wtilde{g}$ nor in $\wtilde{g'}$.
   If it would, then it would appear with the variable $x_{2M}$ in the first bracket in the
   product defining $\wtilde{g}$ and $\wtilde{g'}$, hence there would be a factor
   $$x_{2M-1}^2(x_{2M-1}^n-x_{2M}^n)^2$$
   in the product $\wtilde{g}\cdot\wtilde{g'}$ and then the power of $x_{2M-1}$ would exceed
   $2n$ allowed in $\m$. Hence the variable $x_{2M-1}$ appears in the second bracket in \eqref{eq:g wtilde}.
   Thus we have now
   \begin{equation}\label{eq:g and g' wtilde}
   \begin{array}{rcl}
      \wtilde{g} & = & x_{2M}x_{2M-1}^n\ldots [\ldots x_{2M}][\ldots x_{2M-1}]\\
      \wtilde{g'} & = & x_{2M}x_{2M-1}^n\ldots [\ldots x_{2M}][\ldots x_{2M-1}]
   \end{array}
   \end{equation}
   \item The variable $x_{2M-2}$ in turn has to appear in the first brackets in \eqref{eq:g and g' wtilde}. The argument
   is slightly more involved. In any case there is the factor $(x_{2M-2}^n-x_L^n)$ in $\wtilde{g}$ and $\wtilde{g'}$
   with $L$ either equal to $2M$ or $2M-1$. From these brackets it has to be $x_{2M-2}^n$ which contributes to $\m$
   (otherwise the power at $x_L$ would be too large). So, in any case $x_{2M-2}$ appears with power at least $2n$
   in $\wtilde{g}\cdot \wtilde{g'}$. Since the total power is restricted by $3n$, the only possibility is that
   this variable appears with power $1$ in the products in front of the brackets appearing in \eqref{eq:g and g' wtilde}.
   Hence we have
   \begin{equation}\label{eq:g and g' wtilde 2}
   \begin{array}{rcl}
      \wtilde{g} & = & x_{2M} x_{2M-1}^n x_{2M-2}\ldots [\ldots x_{2M-2} x_{2M}][\ldots x_{2M-1}]\\
      \wtilde{g'} & = & x_{2M} x_{2M-1}^n x_{2M-2}\ldots [\ldots x_{2M-2} x_{2M}][\ldots x_{2M-1}]
   \end{array}
   \end{equation}
   \item Working down, variable by variable, in the same manner, we conclude finally that
   $$\wtilde{g}=\wtilde{g'}=x_2x_4\ldots x_{2M} x_1^nx_3^n\ldots x_{2M-1}^n[x_2x_4\ldots x_{2M}] [x_1x_3\ldots x_{2M-1}]$$
   and the only way to get the monomial $\m$ from the product $\wtilde{g}^2\wtilde{h_{gg}}$ comes in fact
   from the product
   $$x_2^{(2M-2)n+2}x_4^{(2M-4)n+2}\ldots x_{2M}^{2} x_1^{2Mn}x_3^{2(M-1)n}\ldots x_{2M-1}^{2n}\cdot \wtilde{h_{gg}}.$$
   This implies that $\wtilde{h_{gg}}$ contains the monomial
   \begin{equation}\label{eq:monomial p}
      \p=x_2^{n-2}x_4^{n-2}\ldots x_{2M}^{n-2}
   \end{equation}
   with coefficient $1$. But this implies that the coefficient
   of this monomial in $h_{gg}$ is also $1$ (taking modulo $(x_0)$ has no influence on this coefficient).
   \end{itemize}
   The next step is to take \eqref{eq:1} modulo $(x_{2M-1})$. We will denote now the residue class
   of a polynomial $q$ by $\overline{q}$. Thus \eqref{eq:1} becomes
   \begin{equation}\label{eq:1 mod x 2M-1}
      \overline{f}=\sum\limits_{g,g'\in G}\overline{g}\; \overline{g'}\; \overline{h_{gg'}}.
   \end{equation}
   Now we are interested in the monomial
   $$\m'=x_1^{2Nn} x_2^{(2N-1)n}\ldots x_{2M-2}^{3n} x_0^{2n} x_{2N}^n.$$
   Since
   $$\overline{f}=-x_0^nx_1^n\ldots x_{2M-3}^nx_{2M-2}^nx_{2M}^n[x_0x_1\ldots x_{2M-3}x_{2M-2}x_{2M}]$$
   and obviously the monomial $\m'$ comes up in a unique way in the above product,
   its coefficient in $\overline{f}$ is $-1$.

   Running through an analogous procedure as in the reduction modulo $(x_0)$ step,
   we conclude that the monomial $\m'$ appears on the right hand side of \eqref{eq:1 mod x 2M-1}
   only in the square of the generator
   $$\overline{g}=x_2x_4\ldots x_{2M}x_0^n x_1^n x_3^n\ldots x_{2M-3}^n[x_2x_4\ldots x_{2M}][x_0 x_1 x_3\ldots x_{2M-3}]$$
   multiplied by $\overline{h_{gg}}$.
   This shows that the coefficient of the monomial $\p$ defined in \eqref{eq:monomial p} in $\overline{h_{gg}}$ is now $-1$.
   This contradiction shows the assertion
   $$f\not\in I^2.$$

   Now, we study \textbf{the case $N=2M+1$}.
   Assume to the contrary that $f\in I^2$. Then there are
   polynomials $h_{g,g'}$ such that
   \begin{equation}\label{eq:2}
      f=\sum\limits_{g,g'\in G} h_{g,g'} gg'.
   \end{equation}
   Taking \eqref{eq:2} modulo $(x_0)$ we have
   \begin{equation}\label{eq:2 mod x0}
      \wtilde{f}=-x_1^n\ldots x_{2M+1}^n[x_1\ldots x_{2M+1}]=\sum\limits_{g,g'\in G} \wtilde{h_{g,g'}}\cdot \wtilde{g}\cdot \wtilde{g'},
   \end{equation}
   Once again we focus on the coefficient at the monomial
   $$\m=x_1^{(2M+1)n}x_2^{2Mn}\ldots x_{2M}^{2n} x_{2M+1}^n$$
   on both sides of equation \eqref{eq:2 mod x0}. This coefficient is $-1$
   on the left hand side of \eqref{eq:2 mod x0}. It is easy to see that
   there is exactly one way to get this monomial expanding the product defining $\wtilde{f}$.

   By Proposition \ref{prop:generators} a generator $g\in G$ has the form
   $$g=x_{i_0}\ldots x_{i_M}[x_{i_0}\ldots x_{i_M}][x_{j_1}\ldots x_{j_{M+1}}],$$
   with all indices $i_0,\ldots,i_M,j_1,\ldots,j_{M+1}$ mutually distinct.
   If $0\in\left\{i_0,\ldots,i_M\right\}$, then the residue class of $g$ is zero.
   If $0\in\left\{j_1,\ldots,j_{M+1}\right\}$, then the residue class, after possible
   renumbering of indices $\wtilde{g}$, has the form
   $$\wtilde{g}=x_{i_0}\ldots x_{i_M}x_{j_1}^n\ldots x_{j_M}^n[x_{i_0}\ldots x_{i_M}][x_{j_1}\ldots x_{j_M}].$$
   Similarly as in the case of $N=2M$ one can show that there is only one possibility to get the monomial $\m$
   in the right side of the equation \eqref{eq:2 mod x0}.

   This shows that the coefficient of
   $$x_1^{n-2}x_3^{n-2}\ldots x_{2M+1}^{n-2}$$
   in $\wtilde{h_{g,g}}$ (and hence in $h_{g,g}$) is $-1$, where
   \begin{equation}\label{eq: generator g}
      g=x_1x_3\ldots x_{2M+1}[x_1x_3\ldots x_{2M+1}][x_0x_2\ldots x_{2M}].
   \end{equation}

   Finally, we take equation \eqref{eq:2} modulo $(x_{2M})$ and look for the coefficient of
   $$\m'=x_1^{(2M+1)n}x_2^{2Mn}\ldots x_{0}^{2n} x_{2M+1}^n.$$
   Looking at the exponents in $\m'$, we see that there is only one way to obtain this monomial in $\overline{f}$.
   We present here a brief explanation how to produce such a monomial from $\overline{f}$.
   We multiply the following factors
   $$x_i^n[x_0 x_i][x_i x_{i+1}]\ldots [x_i x_{2M+1}] $$
   and take the first element from every bracket except one bracket of the form $[x_0 x_i]$, for which we take the second element. In other words, we proceed as follows
   $$x_i^n[x_0 x_i][x_i x_{i+1}]\ldots [x_i x_{2M+1}] =x_i^n(-x_i^n)[x_i x_{i+1}]\ldots[x_i x_{2M+1}]+\ldots=$$
   $$x_i^n(-x_i^n)(x_i^n)[x_i x_{i+2}]\ldots[x_i x_{2M+1}]+\ldots,$$
   and so on. We do it for all possible $i\in \{1,\ldots, 2M-1\}$ and multiply the results by each other. Finally we multiply all by $-x_0^n[x_0x_{2M+1}]$.
   More precisely we multiply the first element in the bracket and we obtain the monomial $\m'$. Now we calculate the coefficient, which is
   $(-1)^{2M-1}$ from all $[x_0 x_i]$ brackets and one $(-1)$ from $-x_0^n$.

   Summing up we obtain that the coefficient of the monomial
   $$x_1^{n-2}x_3^{n-2}\ldots x_{2M+1}^{n-2}$$
   in $\overline{h_{g,g}}$ (and hence in $h_{g,g}$) where $g$ is as in \eqref{eq: generator g} is $1$, which gives a contradiction.

\endproof

\section{Concluding remarks}
   During preparations of this manuscript we were informed that Ben Drabkin \cite{Dra}
   found another proof of the non-containment Theorem \ref{thm:non-containment}. Since
   his methods are completely different from ours we have decided to include a full
   proof of Theorem \ref{thm:non-containment} also because it reveals particular
   symmetries of the ideals we handle here. We hope to expand this path of thoughts
   in our forthcoming paper \cite{MSS}.

\paragraph{Acknowledgements.}
   We would like to thank Tomasz Szemberg, Marcin Dumnicki and Janusz Gwo\'zdziewicz for helpful suggestions.
   Research of Malara was partially supported by National Science Centre, Poland, grant 2016/21/N/ST1/01491.
   Research of Szpond was partially supported by National Science Centre, Poland, grant 2014/15/B/ST1/02197.


\bigskip \small

\bigskip
   Grzegorz Malara,
   Department of Mathematics, Pedagogical University of Cracow,
   Podchor\c a\.zych 2,
   PL-30-084 Krak\'ow, Poland

\nopagebreak
   \textit{E-mail address:} \texttt{grzegorzmalara@gmail.com}

\bigskip
   Justyna Szpond,
   Department of Mathematics, Pedagogical University of Cracow,
   Podchor\c a\.zych 2,
   PL-30-084 Krak\'ow, Poland

\nopagebreak
   \textit{E-mail address:} \texttt{szpond@up.krakow.pl}

\end{document}